\documentclass[11pt,a4paper,reqno]{amsart}

\pdfoutput=1
\usepackage{amsmath,amsfonts,amssymb,amsthm,color}  
\usepackage{mathrsfs}
\usepackage[utf8]{inputenc}
\usepackage{graphicx}
\usepackage{tikz}
\usetikzlibrary{trees}

\def\R{\mathbb{R}}
\usepackage{aligned-overset}

\usepackage{hyperref}
\hypersetup{
    colorlinks=true,
    linkcolor=blue,
    citecolor=blue,
    filecolor=magenta,      
    urlcolor=cyan,
}

\usepackage[
backend=biber,
style=alphabetic,
sorting=ynt
]{biblatex}
\theoremstyle{plain}
\newtheorem{theorem}{Theorem}[section]

\theoremstyle{definition}

\newtheorem{example}[theorem]{Example}

\numberwithin{equation}{section}
\newtheorem*{theorem*}{Theorem}

\newcommand\numberthis{\addtocounter{equation}{1}\tag{\theequation}}

\title[singular Brascamp-Lieb with dimension 1,2,2, and 1]{On the family of singular Brascamp-Lieb inequalities with dimension datum $(1,2,2,1)$}
\author{Fred Yu-Hsiang Lin}

\begin{document}

\begin{abstract}
	Motivated by the triangular Hilbert transform and the work \cite{kovavc2015dyadic}, we classify a certain family of singular Brascamp-Lieb forms which we associate with the dimension datum $(1,2,2,1)$. We determine the exact range of Lebesgue exponents, for which one has singular Brascamp Lieb inequalities within this family. 
One key observation is a simple proof of a variant of an estimate in \cite{kovavc2015dyadic}.

The remaining observations concern counter examples to boundedness.
  We compare with a counter example showing that the triangular Hilbert form does not satisfy singular Brascamp Lieb bounds with exponents $(\infty, p, p')$.
	\end{abstract}

\maketitle
\section{Introduction}
For $F,G,H\in \mathcal{S}(\mathbb{R}^{2})$, the triangular Hilbert form $\Lambda$ is defined by
\[
\Lambda (F,G,H):=\operatorname{p.v.}\int_{\mathbb{R}^{3}}F(x,y)G(y,z)H(z,x)\frac{dxdydz}{x+y+z}.   \numberthis \label{tht} 
\]
A celebrated open problem in harmonic analysis is whether there exists a constant $C$  such that for all Schwartz functions $F,G,H$, the a priori inequality
\begin{equation}
|\Lambda (F,G,H)| \leq C \|F\|_{L^{p_{1}}(\mathbb{R}^{2})}\|G\|_{L^{p_{2}}(\mathbb{R}^{2})}\|H\|_{L^{p_{3}}(\mathbb{R}^{2})}. \numberthis \label{thtbdd} 
\end{equation}
holds for any exponents $p_1,p_2,p_3$, in particular for the case $p_{1}=p_{2}=p_{3}=3$.
Such bounds are stronger than several important results in harmonic analysis related to Carleson's operator \cite{carleson1966}, \cite{feff1973}, \cite{lt2000} and the bilinear Hilbert transform \cite{lt1997}, \cite{lacey1999calderon} and also the latter's  uniform estimates \cite{thiele2002uniform}, \cite{grafakos2004uniform}, \cite{li2006uniform}, \cite{https://doi.org/10.48550/arxiv.2205.09851}. Inequality (\ref{thtbdd}) falls in the realm of singular Brascamp-Lieb inequalities as in the survey \cite{durcik2021singular}.\\

To make small progress towards the triangular Hilbert form, we discuss a simpler related family of singular Brascamp Lieb forms.  If we consider a function $F$ of the form
\[
 F(x,y):=f(x+\alpha y) \label{specfunc}
 \]
for a one dimensional Schwartz function $f$, then the integral \eqref{tht} is still well defined and becomes
 \[
\widetilde{\Lambda}_{\alpha} (f,G,H):=\operatorname{p.v.}\int_{\mathbb{R}^{3}}f(x+\alpha y)G(y,z)H(z,x)\frac{dxdydz}{x+y+z}.   \numberthis \label{thtsp} 
\]
This is a singular Brascamp-Lieb form with dimensions $1,2,2$ of the functions $f,G,H$  and 
dimension one of the singular kernel.
If $f$ is in $L^\infty(\mathbb{R})$, then $F$ is in $L^\infty(\mathbb{R}^2)$ and singular Brascamp Lieb bounds for \eqref{thtsp} can be 
viewed as special cases of \eqref{thtbdd} with $p_1=\infty$.
For $f$ in $L^p(\R)$  with $p<\infty$,
inequality \eqref{thtsp} is not a special case of inequality \eqref{thtsp}.

More generally, our family of interest of singular Brascamp Lieb forms with dimension datum $(1,2,2,1)$ is a generalization of  \eqref{thtsp} in the spirit of \cite{durcik2021singular} towards  general projections with dimensions $(1,2,2,1)$ in the arguments, that is 
\[
\Lambda_{\Pi}(f,G,H):=\operatorname{p.v.}\int_{\mathbb{R}^{3}}f(\Pi_{1}x)G(\Pi_{2}x)H(\Pi_{3}x)\frac{1}{\Pi_{4}x}dx,
\numberthis \label{genthtsp}
\]
with the projection datum 
\[
\Pi=(\Pi_1,\Pi_2,\Pi_3,\Pi_4),
\numberthis \label{pidatum}
\] where 
\begin{equation}\label{pidim}
\Pi_{1},\Pi_{4}:\mathbb{R}^{3}\rightarrow \mathbb{R},\quad   \Pi_{2},\Pi_{3}:\mathbb{R}^{3}\rightarrow \mathbb{R}^{2}
\end{equation}
are surjective linear maps. To avoid some trivial cases, we assume that
\begin{equation}\label{asspi4}
    \operatorname{Im}(\Pi_4^{\ast
    })\not \subset  \operatorname{Im}(\Pi_j^{*})
\end{equation}
for $j=1,2,3$, where the star denotes the adjoint.

Our first theorem classifies the forms 
\eqref{genthtsp} up to equivalence, 
where $\Lambda_\Pi$ is equivalent
to  $\Lambda_{\Pi'}$, if there exist
\[B\in GL(\mathbb{R}^{3}),\  A_{1},A_{4}\in GL(\mathbb{R}), \  A_{2},A_{3}\in GL(\mathbb{R}^{2}),\] such that for each $1\le j\le 4$
\[\Pi'_j=A_j \Pi_{j} B.\]
Note that then we have
\[
\sup_{f,G,H}\frac{|\Lambda_{\Pi'}(f,G,H)|}{\|f\|_{{p_{1}}}\|G\|_{{p_{2}}}\|H\|_{{p_{3}}}}=\frac{A_{4}^{\frac{1}{p_{1}}+\frac{2}{p_{2}}+\frac{2}{p_{3}}-3}}{\operatorname{detB}\cdot  
\prod_{j=1}^{3}\operatorname{det}(A_{j})^{\frac{1}{p_{j}}}}
\sup_{f,G,H}\frac{|\Lambda_{\Pi}(f,G,H)|}{\|f\|_{{p_{1}}}\|G\|_{{p_{2}}}\|H\|_{{p_{3}}}}.
\]
Hence, it suffices to consider the boundedness of standard forms listed in Theorem \ref{thmclass}.\\

\begin{theorem} \label{thmclass}
    Let $\Pi$ be a datum as in \eqref{pidatum},\eqref{pidim}, and \eqref{asspi4}.
    If $\Lambda_{\Pi}$ is nonzero, it is equivalent to one of 
\begin{align*}
     \Lambda_{(1)}(f,G,H):=&\operatorname{p.v.}\int_{\mathbb{R}^{3}}f(x)G(x,y)H(x,y+t)\frac{1}{t}dtdxdy, \numberthis \label{case1} \\
    \Lambda_{(2)}(f,G,H):=&\operatorname{p.v.}\int_{\mathbb{R}^{3}}f(x)G(x,y)H(x+t,y)\frac{1}{t}dtdxdy \numberthis,  \label{case2}\\ 
    \Lambda_{(3)}(f,G,H):=&\operatorname{p.v.}\int_{\mathbb{R}^{3}}f(x+t)G(x,y)H(x,y+t)\frac{1}{t}dtdxdy, \numberthis  \label{case3}\\
\end{align*}
or it is equivalent to 
\[
\Lambda_{(4,\beta)}(f,G,H):=\operatorname{p.v.}\int_{\mathbb{R}^{3}}f(x+t)G(x,y)H(x+\beta t,y)\frac{1}{t}dtdxdy \numberthis \label{interest}
\]
for some $\beta \in \mathbb{R}$.
Furthermore, any two forms in the above  three discrete cases and in the one parameter family are mutually not equivalent to each other.

\end{theorem}

In particular, Form \eqref{thtsp} is equivalent to
\eqref{case3} if $\alpha=1$ and to \eqref{interest} 
with $\beta=1-\alpha$ otherwise.

Our second goal is to discuss $L^{p}$ bounds of the  forms in Theorem \ref{thmclass}. 
The region of exponents where one may have such bounds is 
\[
    p_{1}=\infty, \:\frac{1}{p_{2}}+\frac{1}{p_{3}}=1, \: 1<p_{2},p_{3}<\infty.
    \numberthis \label{pexp}
\]


Theorem \ref{thmbdd} describes the most interesting case \eqref{interest}, while Theorem \ref{thmother} considers the remaining cases.

\begin{theorem} \label{thmbdd}
    Let  $\beta \neq 0,1$, assume $(p_1,p_2,p_3)$ is in the range \eqref{pexp}. There exists a constant $C$ such that for $f\in \mathcal{S}(\mathbb{R}),\: G,H\in \mathcal{S}(\mathbb{R}^{2})$, the following holds:
    \begin{align*}
|\Lambda_{(4,\beta)}(f,G,H) |
\leq \:C\|f\|_{L^{p_{1}}(\mathbb{R})}\|G\|_{L^{p_{2}}(\mathbb{R}^{2})}\|H\|_{L^{p_{3}}(\mathbb{R}^{2})}.  \numberthis \label{sthtbdd}
    \end{align*}
If $(p_1,p_2,p_3)$ is not in the range \eqref{pexp}, then there is no constant $C$ such that the a priori inequality (\ref{sthtbdd}) holds.

\end{theorem}

The proof of the positive result in this theorem by freezing a variable 
can easily be adapted to the $L^\infty$ case of an estimate in \cite{kovavc2015dyadic} in the dyadic setting.
Somehow this easy case of the family of estimates in \cite{kovavc2015dyadic} appears 
to have been overlooked there.

For the remaining  cases of the classification in Theorem \ref{thmclass}, we summarize the much easier $L^p$ theory in next theorem. 

\begin{theorem}\label{thmother}
    The three variants of Theorem   \ref{thmbdd}   with $\Lambda_{(4,\beta)}$ replaced by $\Lambda_{(1)}$ or $\Lambda_{(2)}$ or $\Lambda_{(4,1)}$ remain true. 
On the other hand, for any $1\leq p_1,p_2,p_3\leq \infty$, there is no constant $C$ such that the two variants of \eqref{sthtbdd} with $\Lambda_{(4,\beta)}$ replaced by $\Lambda_{3}$ or  $\Lambda_{4,0}$ 
hold.
\end{theorem}


We compare the above counter examples with the following, which seems folklore, compared with \cite{kovavc2012boundedness},
but we found useful to do explicitly here.

\begin{theorem}\label{thmcounter}
    Let $\Lambda$ be the triangular Hilbert transform (\ref{tht}). Given exponents $p_{1},p_{2},p_{3}$  satisfying (\ref{pexp}), and $C$ be any constant, there exists $F,G,H\in \mathcal{S}(\mathbb{R}^{2})$ such that
\begin{equation}\label{1.14}
|\Lambda (F,G,H)|\geq C \|F\|_{L^{\infty}}\|G\|_{L^{p_{2}}}\|H\|_{L^{p_{3}}}.
\end{equation}
\end{theorem}



\section*{Acknowledgments}
I am grateful to my advisor Christoph Thiele for posing this problem. I also want to thank Martin Yung-Chang Hsu and Lars Becker for practical discussions on this topic. The author was funded by the DFG under Germany’s Excellence Strategy - EXC-2047/1 - 390685813 and DAAD Graduate School Scholarship Programme - 57572629.

\section{Proof of Theorem \ref{thmclass}}

Let $\Pi$ be a datum as in \eqref{pidatum},\eqref{pidim}, and \eqref{asspi4}. For $1\le j\le 4$, let $V_j$ be the image of $\Pi_j^*$. As $\Pi_j$ is surjective, the
dimension of $V_j$ is one if $j=1,4$, and is two if $j=2,3$. We  do a case distinction according to the relative positions of these subspaces. The following diagram
shows the case distinction.

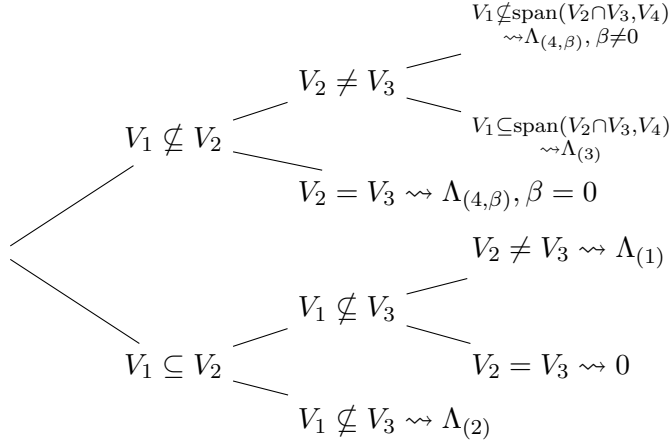
\begin{tikzpicture}[level distance=1.5cm,
  level 1/.style={sibling distance=3cm},
  level 2/.style={sibling distance=1.5cm}]
  \node {}[grow=right,right]
    child {node {$V_{1}\subseteq V_{2}$}
      child {node {$V_{1}\nsubseteq V_{3} \leadsto \Lambda_{(2)}$}}
      child {node {$V_{1}\nsubseteq V_{3}$}
        child {node {$V_{2}=V_{3}\leadsto 0$}}
        child {node {$V_{2}\neq V_{3}\leadsto \Lambda_{(1)}$}}
      }
    }
    child {node {$V_{1}\nsubseteq V_{2}$}
    child {node {$V_{2}=V_{3} \leadsto \Lambda_{(4,\beta)}, \beta =0$ }}
      child {node {$V_{2}\neq V_{3}$}
       child {node {$\substack{V_{1}\subseteq \operatorname{span}(V_{2}\cap V_{3},V_{4})\\ \leadsto \Lambda_{(3)}}$}
       }
       child {node {$\substack{V_{1}\nsubseteq \operatorname{span}(V_{2}\cap V_{3},V_{4})\\ \leadsto \Lambda_{(4,\beta)},\: \beta \neq 0}$}
       }
      }
    };
\end{tikzpicture} 

We first start with the most interesting case.

\textbf{Case 1}: On the one hand $V_{1}\nsubseteq V_{2}$ and on the other hand 
$V_{2}=V_{3}$ or  $V_{1}\nsubseteq \operatorname{span}(V_{2}\cap V_{3}, V_{4})$.\\
As $V_4$ is not contained in $V_2$ by \eqref{asspi4}, we have that $V_2+V_4$ is the full space $\R^3$. Hence we may pick $v_1\in V_2$ and $v_3\in V_4$ such that $v_1+v_3$ spans the one dimensional space $V_1$. The vectors $v_1$ and $v_3$
are linearly independent, because $V_1$ is not in
$V_2$ by the first assumption in Case 1, and not in $V_4$ by \eqref{asspi4}. Now choose a vector $v_2$ in
$V_2\cap V_3$ which is linearly independent of $v_1$.
This is possible if $V_2=V_3$, because then $V_{2}\cap V_{3}$ is two dimensional.
It is also possible if $V_{1}\nsubseteq \operatorname{span}(V_{2}\cap V_{3}, V_{4})$.
Namely, let $v_2$ be any nonzero vector in $V_2\cap V_3$ and assume to get a contradiction that $v_1$ is  a multiple
of $v_2$. Then $v_1+v_3$ is in $\operatorname{span}(V_{2}\cap V_{3}, V_{4})$. This contradicts that $v_1+v_3$ spans $V_1$. Hence we have seen that under the assumption of Case 1, we can choose $v_{2}$ as above.

As $v_1$ and $v_3$ are linearly independent and $v_3$ is not in $V_3$, there is a $\beta\in \R$ such that $v_1+\beta v_3\in V_3$. Since $V_{2}+V_{4}=\mathbb{R}^{3}$, $\operatorname{span}\{v_{1},v_{2}\}=V_{2}$, and $\operatorname{span}\{v_{3}\}=V_{4}$, we have found a basis $\{v_{1},v_{2},v_{3}\}$ for $\mathbb{R}^{3}$.

We choose $B$ so that $B^*$ maps 
$v_1,v_2,v_3$ to the standard unit vectors $e_1,e_2,e_3$. We choose $A_j$ so that $A_j^*$ maps the standard unit vectors of $\R^1$ or $\R^2$ to preimages under $\Pi_j^*$
of the spanning vectors of $V_j$ expressed as above in terms of  $v_1$, $v_2$, or $v_3$. 
This allows to directly write down the matrix for $B^*\Pi_j^*A_j^*$ in the standard basis. 
The transposes of these matrices are as follows.

\[A_{1}\Pi_{1}B=
\begin{pmatrix}
   1&0&1
\end{pmatrix},
A_{2}\Pi_{2}B=
\begin{pmatrix}
  1&0&0 \\0&1&0
\end{pmatrix},
\]
\[
A_{3}\Pi_{3}B=
\begin{pmatrix}
   1&0& \beta \\0&1&0
\end{pmatrix},
A_{4}\Pi_{4}B=
\begin{pmatrix}
  0&0&1
\end{pmatrix}.
\]
Hence 
\[
\operatorname{p.v.}\int_{\mathbb{R}^{3}}f(
A_1\Pi_{1}BX)G(A_2\Pi_{2}BX)H(A_3\Pi_{3}BX)\frac{1}{A_4\Pi_{4}BX}dX
\numberthis \label{bigx}
\]

\[
=\operatorname{p.v.}\int_{\mathbb{R}^{3}}f(x+t)G(x,y)H(x+\beta t,y)\frac{1}{t}dtdxdy  . 
\]
which is $\Lambda_{(4,\beta)}$.\\

\textbf{Case 2}: $V_{1}\nsubseteq V_{2}$, $V_{2}\neq V_{3}$, and $V_{1}\subseteq \operatorname{span}(V_{2}\cap V_{3}, V_{4})$.\\

By the last assumption of Case 2, pick $v_{1}\in V_{2}\cap V_{3}$ and $v_{3}\in V_{4}$ such that $v_{1}+v_{3}$ spans $V_{1}$. The vectors $v_{1},v_{3}$ are linearly independent, otherwise one of them is a nonzero vector in $V_{1}$ which is impossible by  \eqref{asspi4} and the first assumption of Case 2. Choose $v_{2}$ to be a vector in $V_{2}$ which is linearly independent of $v_{1}$ such that $v_{2}+v_{3}$ is a nonzero vector in $ V_{3}$. This is possible because $V_{2}\neq V_{3}$ by the second assumption of Case 2 and $V_{4}\cap V_{3}= \{0\}$ by \eqref{asspi4}. Since $V_{1}\nsubseteq V_{2}$, we have $v_{2}\notin \operatorname{span}\{v_{1},v_{3}\}$. Hence $\{v_{1},v_{2},v_{3}\}$ forms a basis of $\mathbb{R}^{3}$. 

Choosing $A_j$ and $B$ similarly as above, we have
\[A_{1}\Pi_{1}B=
\begin{pmatrix}
    1&0&1
\end{pmatrix},
A_{2}\Pi_{2}B=
\begin{pmatrix}
   1&0&0 \\ 0&1&0
\end{pmatrix},
\]
\[
A_{3}\Pi_{3}B=
\begin{pmatrix}
    1&0&0 \\ 0&1&1
\end{pmatrix},
A_{4}\Pi_{4}B=
\begin{pmatrix}
   0&0&1
\end{pmatrix}.
\]
Hence \eqref{bigx} is equal to
\[\operatorname{p.v.}\int_{\mathbb{R}^{3}}f(x+t)G(x,y)H(x,y+t)\frac{1}{t}dtdxdy  . 
\]
which is $\Lambda_{(3)}$.\\

\textbf{Case 3}: $V_{1}\subseteq V_{2}$ and $V_{1}\nsubseteq V_{3}$.\\
Let $v_{2}$ be a vector in $V_{2}\cap V_{3}$. Choose $v_{1}\in V_{1}$ and $v_{3}\in V_{4}$ such that $v_{1}+v_{3}\in V_{3}$. This is possible because $V_{1}\nsubseteq V_{3}$ and $V_{4}\nsubseteq V_{3}$. Since $V_{4}\nsubseteq V_{2}$ and $V_{1}\subseteq V_{2}$, we have $V_{2}\cap V_{3} \nsubseteq V_{1}+V_{4}$. This shows that $\{v_{1},v_{2},v_{3}\}$ is a basis of $\mathbb{R}^{3}$.
Choosing $A_j$ and $B$ similarly as above, we have
\[
A_{1}\Pi_{1}B=
\begin{pmatrix}
   1&0&0
\end{pmatrix},
A_{2}\Pi_{2}B=
\begin{pmatrix}
   1&0&0 \\0&1&0
\end{pmatrix},
\]
\[
A_{3}\Pi_{3}B=
\begin{pmatrix}
   1&0&1 \\0&1&0
\end{pmatrix},
A_{4}\Pi_{4}B=
\begin{pmatrix}
   0&0&1
\end{pmatrix}.
\]
Hence \eqref{bigx} is equal to
\[
\operatorname{p.v.}\int_{\mathbb{R}^{3}}f(x)G(x,y)H(x+t,y)\frac{1}{t}dtdxdy  . 
\]
which is $\Lambda_{(2)}$.\\

\textbf{Case 4}:  $V_{1}\subseteq V_{2}$, $V_{1}\subseteq V_{3}$, and $V_2\neq V_3$.\\
Let $v_{1}$ be a vector in $V_{1}$ and $v_{3}$ be a vector in $V_{4}$. Choose another vector $v_{2}\in V_{2}$ which is linearly independent of $v_{1}$ such that $v_{2}+v_{3}\in V_{3}$. This is possible because $V_{2}\neq V_{3}$ and $V_{4}\nsubseteq V_{3}$. Since $V_{2}+V_{4}=\mathbb{R}^{3}$, we have $\{v_{1},v_{2},v_{3}\}$ forms a basis of $\mathbb{R}^{3}$.

Choosing $A_j$ and $B$ suitably as above, we have
\[A_{1}\Pi_{1}B=
\begin{pmatrix}
    1&0&0
\end{pmatrix},
A_{2}\Pi_{2}B=
\begin{pmatrix}
   1&0&0 \\0&1&0
\end{pmatrix},
\]
\[
A_{3}\Pi_{3}B=
\begin{pmatrix}
    1&0&0 \\0&1&1
\end{pmatrix},
A_{4}\Pi_{4}B=
\begin{pmatrix}
    0&0&1
\end{pmatrix}.
\]
Similarly as above, we obtain for \eqref{bigx} 
\[
\operatorname{p.v.}\int_{\mathbb{R}^{3}}f(x)G(x,y)H(x,y+t)\frac{1}{t}dtdxdy  . 
\]
which is $\Lambda_{(1)}$.

Notice that beside the above cases, there's a trivial case, $V_{1}\subseteq V_{2}$ and $V_2=V_3$. We then have the following data.

\[A_1 \Pi_1 B=
\begin{pmatrix}
    1 & 0 & 0
\end{pmatrix},\quad 
A_2 \Pi_2 B=
\begin{pmatrix}
    1& 0 &0 \\ 0& 1 & 0
\end{pmatrix},\]
\[A_{3}\Pi_{3}B=
\begin{pmatrix}
   1&0&0 \\0&1&0
\end{pmatrix},\quad
A_{4}\Pi_{4}B=
\begin{pmatrix}
   0&0&1
\end{pmatrix}.
\]

and the corresponding trilinear form
\[
\operatorname{p.v.}\int_{\mathbb{R}^{3}}f(x)G(x,y)H(x,y)\frac{1}{t}dtdxdy  . 
\]

This integral is zero since the principal value is interpreted as limit as $\epsilon \to 0$ 
of the truncation of the integral
to $t\in [-\epsilon^{-1}, \epsilon^{-1}]\setminus [-\epsilon, \epsilon]$.

Since the inclusion relation of subspaces maintain the same after basis change, $\Lambda_{(1)}$, $\Lambda_{(2)}$, $\Lambda_{(3)}$, $\Lambda_{(4,0)}$, and $\Lambda_{(4,1)}$ are mutually not equivalent to each other and all not equivalent to $\Lambda_{(4,\beta)}$ for $\beta \neq 0,1$. In the following, we define a quantity in projective geometry to distinguish $\Lambda_{(4,\beta)}$ for different $\beta$. Let $V_{1}'$ be the space spanned by $V_{1}$ and $V_{2}\cap V_{3}$. Let $V_{4}'$ be the space spanned by $V_{4}$ and $V_{2}\cap V_{3}$. We define the cross ratio of these four planes as follows. Take an arbitrary line $L$ in $\mathbb{R}^{3}$ not intersecting $V_{2}\cap V_{3}$. Let $x_{1}, x_{2}, x_{3},$ and $x_{4}$ be the intersection of $L$ with these four planes then project to $\mathbb{R}$ respectively. Then the cross ratio
\[
\frac{(x_{1}-x_{4})(x_{3}-x_{2})}{(x_{1}-x_{2})(x_{3}-x_{4})}
\]
is independent of the choice of $L$ and is an invariant under basis change. We calculate the cross ratio for $\Lambda_{4,\beta}$ with $\beta \neq 0,1$. Consider the line $L:x+z=1,\: y=0$. The intersection of $L$ with these four planes are $(\frac{1}{2},0,\frac{1}{2})$, $(1,0,0)$, $(\frac{1}{1+\beta},0,\frac{\beta}{1+\beta})$, $(0,0,1)$ respectively. To calculate the cross ratio of these four points, it suffices to focus on its $x$ variable.
\[
\frac{(\frac{1}{1+\beta}-1)\cdot \frac{1}{2}}{\frac{1}{1+\beta}\cdot (-\frac{1}{2})}=\beta.
\]
With this invariant, we then finish all the classification in the last case.

\begin{example}\label{example}
    We will write (\ref{thtsp}) in the standard form.

For $\alpha\neq 1$, (\ref{thtsp}) is equivalent to $\Pi_{(4,1-\alpha)}$. For  $\alpha=1$ (\ref{thtsp}) is equivalent to $\Pi_{(3)}$. Moreover, we can explicitly show $A_{1},A_{2},A_{3},A_{4},$ and  $B$. Let $\Pi$ be the datum of (\ref{thtsp}).

For $\alpha \neq 1$, take
\[
A_{1}=1,\quad A_{2}=
\begin{pmatrix}
    \alpha -1 &-1\\
    0&1
\end{pmatrix},\quad 
 A_{3}=
\begin{pmatrix}
    -\alpha  &-\alpha +1\\
    1&0
\end{pmatrix},\quad
A_{4}=1,\numberthis \label{aneq1m}
\]
and 
\[
B=
\begin{pmatrix}
    \frac{-1}{\alpha -1}&\frac{-\alpha}{\alpha -1}&1\\
    \frac{1}{\alpha -1}&\frac{1}{\alpha -1}&0\\
    0&1&0
\end{pmatrix}. \numberthis \label{Baneq1m}
\]
Then
\[A_{1}\Pi_{1}B=
\begin{pmatrix}
    1& 0 &1
\end{pmatrix},
A_{2}\Pi_{2}B=
\begin{pmatrix}
    1& 0&0\\
    0&1&0
\end{pmatrix},
\]
\[
A_{3}\Pi_{3}B=
\begin{pmatrix}
     1& 0& 1-\alpha\\
    0&1&0
\end{pmatrix},
A_{4}\Pi_{4}B=
\begin{pmatrix}
    0&0&1
\end{pmatrix}.
\]
For $\alpha =1$, take
\[
A_{1}=1,\quad
A_{2}=
\begin{pmatrix}
    0&1\\
    1&0
\end{pmatrix},\quad
A_{3}=
\begin{pmatrix}
    1&0\\
    -1&-1
\end{pmatrix},\quad
A_{4}=1, \numberthis \label{aeq1m}
\]
and
\[
B=
\begin{pmatrix}
    -1&-1&1\\
    0&1&0\\
    1&0&0
\end{pmatrix}.    \numberthis \label{Baeq1m}
\]
Then
\[A_{1}\Pi_{1}B=
\begin{pmatrix}
    1& 0 &1
\end{pmatrix},
A_{2}\Pi_{2}B=
\begin{pmatrix}
    1& 0&0\\
    0&1&0
\end{pmatrix},
\]
\[
A_{3}\Pi_{3}B=
\begin{pmatrix}
     1& 0& 0\\
    0&1&1
\end{pmatrix},
A_{4}\Pi_{4}B=
\begin{pmatrix}
    0&0&1
\end{pmatrix}.
\]
\end{example}

\section{Proof of Theorem \ref{thmbdd}} 

 Let $\beta \in \mathbb{R}$, assume $(p_1,p_2,p_3)$ is in the range \eqref{pexp}.
The letter $C$ will denote a sufficiently large positive number 
that may be implicitly re-adjusted from inequality to inequality 
and that may depend on $\beta$ and $p_1,p_2,p_3$.
We write $A\lesssim B$ if $A\leq CB$ for such number $C$. We write $A\sim B$ if both $A\lesssim B$  and $B\lesssim A$. We will adopt this convention in the rest of this paper.

First consider exponents $(p_{1},p_{2},p_{3})$ in the range (\ref{pexp}). 
Assuming momentarily  that we can pass the $\operatorname{p.v}$ inside, then freezing the $y$ variable, we may identify the form as the trilinear form associated with the bilinear Hilbert transform and obtain the desired estimate (\ref{sthtbdd}) as follows:

\begin{align*}
    &\left|\operatorname{p.v.}\int_{\mathbb{R}^{3}}f(x+t)G(x,y)H(x+\beta t,y)\frac{1}{t}dtdxdy \right|\\
=&\left|\int_{\mathbb{R}}\:\:\operatorname{p.v.}\int_{\mathbb{R}^{2}}f(x+t)G(x,y)H(x+\beta t,y)\frac{1}{t}dtdxdy\right|\\
    \lesssim &\int_{\mathbb{R}}\|f\|_{L^{\infty}}\|G_{y}\|_{L^{p_{2}}}\|H_{y}\|_{L^{p_{3}}}dy\\
    \lesssim 
    &\|f\|_{L^{\infty}}\|G\|_{L^{p_{2}}}\|H\|_{L^{p_{3}}}.
\end{align*}

To show that we can pass the $\operatorname{p.v.}$ inside, by dominant convergence theorem, it suffices to show that
$M\in L^1(\R)$, where
\[
M(y):=\operatorname{sup}_{\varepsilon >0}\left|\int_{|t|>\varepsilon}f(x+t)G(x,y)H(x+\beta t,y)\frac{1}{t}dtdx\right| .\numberthis \label{pvpast1}
\]
We set \ $G_{y}(x):=G(x,y),\: H_{y}(x):=H(x,y)$. We split the integrand into $|t|>1$ and $|t|<1$. By triangle inequality  we have
\[
\begin{aligned}
|M(y)|&\leq \int_{|t|\leq 1}\left|\frac{f(x+t)-f(x)}{t}G_{y}(x)H_{y}(x+\beta t)\right|dtdx\\
&+\int_{|t|\leq 1}\left|f(x)G_{y}(x)\frac{H_{y}(x+\beta t)-H_{y}(x)}{t} \right|dtdx\\
&+\int_{|t|\leq 1}\left|f(x)G_{y}(x)H_{y}(x)\right|dtdx\\
&+\int_{1<|t|}\left| f(x+t)G_{y}(x)H_{y}(x+\beta t)\frac{1}{t}\right|dtdx.
\end{aligned}
\]
Then by mean value theorem,
\[
\begin{aligned}
    &\leq \int_{|t|\leq 1}\left|\|f'\|_{L^{\infty}}G_{y}(x)H_{y}(x+\beta t)\right|dtdx\\
&+\int_{|t|\leq 1}\left|f(x)G_{y}(x)\|H_{y}'\|_{L^{\infty}} \right|dtdx\\
&+\int_{|t|\leq 1}\left|f(x)G_{y}(x)H_{y}(x) \right|dtdx\\
&+\int_{1<|t|}\left| f(x+t)G_{y}(x)H_{y}(x+\beta t)\right|dtdx.
\end{aligned}
\]

Since $f$, $G$, $H$ are all Schwartz function, these four terms are all $L^{1}$ integrable. 
This completes the proof of estimate (\ref{sthtbdd}) for 
$(p_1,p_2,p_3)$ in the range \eqref{pexp}.

Now assume $(p_{1},p_{2},p_{3})$ do not satisfy (\ref{pexp}), we will show that  the a priori inequality (\ref{sthtbdd}) does not hold. In this proof, we will use $H$ to denote the Hilbert transform. To distinguish from the Hilbert transform, we take $E$ to denote the third function in our trilinear form. We will prove by contradiction.

\begin{align*}
    \Lambda_{\Pi}(T_{k}g,G,T_{\beta k}^{(1)}E)&=\operatorname{p.v.}\int_{\mathbb{R}^{3}}f(x-k+t)G(x,y)E(x-\beta k+\beta t,y)\frac{1}{t}dtdxdy\\
    &=\operatorname{p.v.}\int_{\mathbb{R}^{3}}f(x+t)G(x,y)E(x+\beta t,y)\frac{1}{t+k}dtdxdy.
\end{align*}
Now fix a $m\in L^{1}$ with $\operatorname{supp}\widehat{m}\subseteq (0,\infty)$. 
\begin{align*}
    &|\int_{\mathbb{R}}(Hm)(-k)\Lambda_{\Pi}(T_{k}g,G,T_{\beta k}^{(1)}E)dk|\\
    \lesssim &\int_{\mathbb{R}}|Hm(-k)|\cdot |\Lambda_{\Pi}(T_{k}g,G,T_{\beta k}^{(1)}E)|dk\\
    \lesssim &\int_{\mathbb{R}}|Hm(-k)|\cdot \|T_{k}f\|_{L^{p_{1}}}\|G\|_{L^{p_{2}}}\|T_{\beta k}^{(1)}E\|_{L^{p_{3}}}dk\\
    \lesssim &\int_{\mathbb{R}}|Hm(-k)|\cdot \|f\|_{L^{p_{1}}}\|G\|_{L^{p_{2}}}\|E\|_{L^{p_{3}}}dk.
\end{align*}
By $Hm=-im$, where $H$ is the Hilbert transform, 
\[
\left|\int_{\mathbb{R}}(Hm)(-k)\Lambda_{\Pi}(T_{k}f,G,T_{\beta k}^{(1)}E)dk\right| \lesssim \|f\|_{L^{p_{1}}}\|G\|_{L^{p_{2}}}\|E\|_{L^{p_{3}}} \|m\|_{L^{1}}. \numberthis \label{quabdd}
\]
On the other hand, by $H^{2}m=-m$, i.e.
\[
m(x)=-\operatorname{p.v.}\int_{\mathbb{R}}(Hm)(-y)\frac{1}{x+y}dy.
\]
We obtain
\begin{align*}
    &\left|\int_{\mathbb{R}}(Hm)(-k)\Lambda_{\Pi}(T_{k}g,G,T_{\beta k}^{(1)}E)dk\right|\\
    =&\left|\int_{\mathbb{R}^{3}}f(x+t)G(x,y)E(x+\beta t,y)m(t)dtdxdy\right|. \numberthis 
\end{align*}
It's temptating to check the Brascamp-Lieb conditions in \cite{bennett2008brascamp} at this point.\\
For scaling condition, 
\[
\frac{1}{p_{1}}+\frac{2}{p_{2}}+\frac{2}{p_{3}}+1=3.  \numberthis \label{scaleond}
\]
As for dimension condition, we may  take the subspace\\ $\operatorname{span}\{(0,1,0)\}$,
\[
1\leq \frac{0}{p_{1}}+\frac{1}{p_{2}}+\frac{1}{p_{3}}+0. \numberthis \label{dimcond}
\]
However, in this case, $m$ is not an arbitrary function. We cannot take $m$ as Gaussian, hence we need to mimic the proof in \cite{bennett2008brascamp} and modify it slightly. Take $f(x)=e^{-\pi x^{2}}$, $G(x,y)=e^{-\pi (x^{2}+\varepsilon y^{2})}$, $E(x,y)=e^{-\pi(x^{2}+\varepsilon y^{2})}$, $\widehat{m}$ is a nonnegative bump function support in $(0,\infty)$. Then
\[
\|f\|_{L^{p_{1}}}\|G\|_{L^{p_{2}}}\|E\|_{L^{p_{3}}}\|m\|_{L^{1}}\sim \varepsilon^{-\frac{1}{2}(\frac{1}{p_{2}}+\frac{1}{p_{3}})}. \numberthis \label{1over2ep}
\]
On the other hand
\[
\int_{\mathbb{R}^{3}}f(x+t)G(x,y)E(x+\beta t,y)m(t)dtdxdy
\]
\[
=\int_{\mathbb{R}^{3}}e^{-\pi (2\varepsilon y^{2})}\cdot e^{-\pi (3[x+\frac{1}{3}(1+\beta t)]^{2})}\cdot e^{-\pi (\frac{1}{3}(2\beta^{2}+\beta+2)t^{2})}m(t)dydxdt. \numberthis \label{eplong}
\]
First integrate in $y$, we may obtain a factor $\varepsilon^{-\frac{1}{2}}$. Second, integrate in $x$, we may get a constant. Hence, we may simplify (\ref{eplong}) into
\[
\sim \varepsilon^{-\frac{1}{2}} \int_{\mathbb{R}^{3}} e^{-\pi (\frac{1}{3}(2\beta^{2}+\beta+2)t^{2})}m(t)dt. \numberthis \label{eplong2}
\]
By Plancherel identity, we may further reduce (\ref{eplong2}) to
\[
\varepsilon^{-\frac{1}{2}}(\frac{1}{3}(2\beta^{2}+\beta+2))^{-\frac{1}{2}} \int_{\mathbb{R}^{3}} e^{-\pi ((3\beta^{2}+\beta+2)^{-1}t^{2})}\widehat{m}(t)dt. \numberthis \label{eplong3}
\]
The integral in (\ref{eplong3}) is a positive constant away from $0$ and $2\beta^{2}+\beta+2$ is positive for all $\beta \in \mathbb{R}$. Taking $\varepsilon \rightarrow 0$ and compare (\ref{1over2ep}) and (\ref{eplong3}), we may have the exponents $p_{1},p_{2},p_{3}$ in the quadrilinear form (\ref{quabdd}) is impossible to hold for the range other than (\ref{pexp}), which is a contradiction. Hence (\ref{sthtbdd}) cannot hold for exponents $(p_{1},p_{2},p_{3})$ other than (\ref{pexp}).\\
This completes  the proof of Theorem \ref{thmbdd}.

\section{Proof of Theorem \ref{thmother}}
Now we discuss the $L^{p}$ bounded for the form (\ref{genthtsp}). 
Since equivalent forms share the same $L^{p}$ boundedness property,  suffice to consider the boundedness of standard forms listed in Theorem \ref{thmclass}.\\

First, we deal with the endpoint range (\ref{pexp}). The cases (\ref{case3}) and (\ref{interest}) with $\beta=0$ are not bounded. On the other hand, the cases (\ref{case1}), (\ref{case2}), and (\ref{interest}) with $\beta =1$ are bounded.\\

For (\ref{case3}) and (\ref{interest}) with $\beta=0$, we prove the unboundedness for the case (\ref{case3}), the others are similar. Suppose we have the bound, take 
\[
G(x,y)=\operatorname{sgn}g(x)|g(x)|^{\frac{1}{2}}(D_{N}^{\infty}\varphi )(y),\: H(x,y)=|g(x)|^{\frac{1}{2}}(D_{N}^{\infty}\varphi )(y).
\]
Then 
\begin{align*}
    &\left|\operatorname{p.v.} \int_{\mathbb{R}^{3}}f(x+t)g(x)(D_{N}^{\infty}\varphi )(y)(D_{N}^{\infty}\varphi )(y+t)\frac{1}{t}dxdy\right|\\
    \lesssim  &\|f\|_{L^{\infty}}\cdot (\|g\|_{L^{1}}^{\frac{1}{2}}N^{\frac{1}{2}})\cdot (\|g\|_{L^{1}}^{\frac{1}{2}}N^{\frac{1}{2}})\\
    =&N \|f\|_{L^{\infty}}\cdot \|g\|_{L^{1}}. \numberthis
\end{align*}
Notice that taking $N\rightarrow \infty$ and integrate over $y$, the quantity
\[
\frac{1}{N}(D_{N}^{\infty}\varphi )(y)(D_{N}^{\infty}\varphi )(y+t)
\]
will tend to a constant. This implies Hilbert transform is bounded at $L^{\infty}$, a contradiction.\\

For the cases (\ref{case1}), (\ref{case2}), and (\ref{interest}) with $\beta =1$, we prove the boundedness for the case (\ref{interest}) with $\beta =1$, the proof for other cases are similar.
\begin{align*}
    &\left|\operatorname{p.v.}\int_{\mathbb{R}^{3}}f(x+t)G(x,y)H(x+t,y)\frac{1}{t}dtdxdy\right|\\
    \lesssim &\int_{\mathbb{R}}\|fH_{y}\|_{L^{p_{2}}}\|G_{y}\|_{L^{p_{3}}}dy\\
    \lesssim &\|f\|_{L^{\infty}}\int_{\mathbb{R}}\|H_{y}\|_{L^{p_{2}}}\|G_{y}\|_{L^{p_{3}}}dy\\
    \lesssim &\|f\|_{L^{\infty}}\|G\|_{L^{p_{2}}}\|H\|_{L^{p_{3}}}. \numberthis
\end{align*}
First freeze $y$, pair $f$ and $G_{y}$ together and use the estimate of Hilbert transform. Then use H\"older inequality twice to get the desired estimate.\\

Second, for the impossibility of $L^{p}$ bound of all other cases for exponents $(p_{1},p_{2},p_{3})$ other than the range (\ref{pexp}), the proof is similar to the case in Theorem \ref{thmbdd}.

\section{Proof of Theorem \ref{thmcounter}}

Let exponents $p_{1},p_{2},p_{3}$  be given satisfying (\ref{pexp}), and let $C$ be any constant.
Shearing some functions, we rewrite the triangular Hilbert transform
\begin{equation}\label{5.1}
\begin{aligned}
    \Lambda (F,G,H)&=\operatorname{p.v.}\int_{\mathbb{R}^{3}}F(t-y-z,y)G(y,z)H(z,t-y-z)\frac{1}{t}dtdydz\\
    &=\operatorname{p.v.}\int_{\mathbb{R}^{3}}\widetilde{F}(z-t,y)G(y,z)\widetilde{H}(z,y-t)\frac{1}{t}dtdydz, 
\end{aligned}
\end{equation}
where $\widetilde{F}(x,y)=F(-x-y,y)$, $\widetilde{H}(x,y)=H(x,-y-z)$.\\
For a parameter $N> 2e^{2C}$, let  
\[G_N(y,z)= 1_{[0,N]}(y)1_{[0,1]}(z)\]
\[
\widetilde{H}_{N}(z,y):=1_{[0,N]}(y)1_{[0,1]}(z)
\]
\[
\widetilde{F}(x,y):= 1_{[-\infty,-1)}(x)
\]
Then the integrand of  \eqref{5.1} is non-negative, and we may estimate  \eqref{5.1} from below as
\[
\Lambda (F,G,H)\geq \int_{0}^{1}\int_{\mathbb{R}}\int_{1\leq t}G(y,z)\widetilde{H}(z,y-t)\frac{1}{t}dtdydz. \numberthis
\]
\begin{equation}\label{5.3}
\ge 
\int_0^{\frac N2}\int_{1}^{\frac N2}\frac{1}{t}dtdy=\frac N2 \log(\frac N2)>NC.
\end{equation}

From penultimate to ultimate line,  the we integrated $z$ from $0$ to $1$ and used that if $0<y,t<\frac{N}{2}$, then
both $y$ and $y+t$ are in $[0,N]$.
On the other hand, as shearing leaves the $L^p$ norm invariant, 
\begin{equation}\label{5.4}
\|F\|_\infty \|G\|_{p_2}\|H\|_{p_3}=
\|\tilde{F}\|_\infty \|G\|_{p_2}\|\tilde{H}\|_{p_3}
=N^{\frac 1{p_2}} N^{\frac 1{p_3}}=N.
\end{equation}
This together with \eqref{5.3} shows \eqref{1.14} and completes the proof of Theorem \eqref{thmcounter}.

\printbibliography[heading=bibintoc,title={References}]
\end{document}